\documentclass[12pt]{article}
\usepackage{amssymb,amsmath}
\usepackage{cases}
\usepackage{amsfonts}
\usepackage{color,xcolor}
\usepackage[left=2.0cm,right=2.0cm,top=2.0cm,bottom=2.0cm]{geometry}
\usepackage[colorlinks,citecolor=blue,urlcolor=blue]{hyperref}

\newtheorem{theorem}{Theorem}[section]

\newtheorem{lemma}{Lemma}[section]

\newtheorem{definition}{Definition}[section]
\newtheorem{remark}{Remark}[section]

\newcommand{\bal}{\begin{align}}
\newcommand{\bbal}{\begin{align*}}
\newcommand{\beq}{\begin{equation}}
\newcommand{\eeq}{\end{equation}}
\newcommand{\bca}{\begin{cases}}
\newcommand{\eca}{\end{cases}}

\begin{document}

\title{Ill-posedness for a two-component Novikov system  in Besov space}

\author{Xing Wu$^{1, }$\thanks{Corresponding author. ny2008wx@163.com (Xing Wu)}\;, Min Li$^2$\\
\small \it$^1$College of Information and Management Science,
Henan Agricultural University,\\
\small  Zhengzhou, Henan, 450002, China\\
\small \it $^2$Department of mathematics, Jiangxi University of Finance and Economics,\\
\small Nanchang, Jiangxi, 330032, China}

\date{}

\maketitle\noindent{\hrulefill}

{\bf Abstract:} In this paper, we consider the Cauchy problem for a two-component Novikov system on the line. By specially constructed initial data $(\rho_0, u_0)$ in $B_{p, \infty}^{s-1}(\mathbb{R})\times B_{p, \infty}^s(\mathbb{R})$ with $s>\max\{2+\frac{1}{p}, \frac{5}{2}\}$ and $1\leq p \leq \infty$, we show that  any energy bounded solution starting from $(\rho_0, u_0)$ does not converge back to $(\rho_0, u_0)$ in the metric of $B_{p, \infty}^{s-1}(\mathbb{R})\times B_{p, \infty}^s(\mathbb{R})$ as time goes to zero, thus results in discontinuity of the data-to-solution map and ill-posedness.


\vskip0mm\noindent{\hrulefill}

\section{Introduction}\label{sec1}
The two-component Novikov system takes the form
\begin{eqnarray}\label{eq1}
        \left\{\begin{array}{ll}
         \rho_t=\rho_xu^2+\rho uu_x,\;\; t>0, \;\; x\in \mathbb{R},\\
          m_t=3u_xum+u^2m_x-\rho(u\rho)_x, \;\; t>0, \;\; x\in \mathbb{R},\\
          m=u-u_{xx},\\
         \rho(0, x)=\rho_0, u(0, x)=u_0, \end{array}\right.
        \end{eqnarray}
which was recently proposed by Popowicz  \cite{Popowicz 2015} as  the two-component generalization of the Novikov equation and can be rewritten in the Hamiltonian form (see \cite{Popowicz 2015} for details).

We recall that system (\ref{eq1}) is well-posed in the sense of Hadamard (see \cite{Bahouri 2011}) in Besov spaces $B_{p, r}^{s-1} \times B_{p, r}^s$ with $s>\max\{1+\frac{1}{p}, \frac{3}{2}\}$, $1\leq p\leq \infty$, $1\leq r < \infty$ as well as in the critical Besov space $B_{2, 1}^{\frac{1}{2}}\times B_{2, 1}^{\frac{3}{2}}$ (see \cite{Luo 2015}). More precisely, if $(\rho_0, u_0)$ belongs to $B_{p, r}^{s-1} \times B_{p, r}^s$ on the line, then there exists $T=T(\|\rho_0\|_{B_{p, r}^{s-1} }, \|u_0\|_{B_{p, r}^s})$, and a unique solution $(\rho, u) \in L^\infty(0, T; B_{p, r}^{s-1} \times B_{p, r}^s)$ satisfying system (\ref{eq1}) with continuous dependence on initial data. Note that if we want to include the endpoint case $r=\infty$, then we need to weaken the notion of well-posedness since the continuity of data-to-solution map in $B_{p, \infty}^{s-1}\times B_{p, \infty}^s$ is still unkown.

If we set $\rho=0$, then system (\ref{eq1}) becomes the famous Novikov equation
\begin{equation}\label{eq2}
 m_t=3u_xum+u^2m_x, \qquad m=u-u_{xx},
\end{equation}
which was derived by Novikov \cite{Novikov 2009} as a new integrable equation with cubic nonlinearities. It is shown the Novikov equation is  integrable in the sense of having a Lax pair in matrix form, a bi-Hamiltonian structure as well as possessing infinitely many conserved quantities \cite{Home 2008}. Furthermore, it admits peakon solutions given by the formula $u(x, t)=\pm \sqrt{c}e^{|x-ct|}$ ($c$=const. is the wave speed) and also multipeakon traveling wave solutions.

The local well-posedness for Novikov equation was initially established in the Sobolev spaces $H^s$ with $s>\frac{3}{2}$ on both the line and the circle by Himonas and Holliman \cite{Himonas 2012}, and then by Ni and Zhou \cite{Ni 2011} in the critical Besov space $B_{2, 1}^{\frac{3}{2}}(\mathbb{R})$. Later, the well-posed space was extended to a larger class of Besov spaces $B_{p, r}^s(\mathbb{R})$, $s>\max\{1+\frac{1}{p}, \frac{3}{2}\}$, $1\leq p\leq \infty$, $1\leq r< \infty$ by Yan, Li and Zhang \cite{Yan 2012}. Moreover, by using the peakon traveling wave solution they showed that the local well-posedness fails in $B_{2, \infty}^{\frac{3}{2}}(\mathbb{R})$ due to failure of continuity.
Several years later, by constructing a 2-peakon
solution with an asymmetric antipeakon-peakon,  Himonas, Holliman and Kenig \cite{Himonas 2018} proved ill-posedness for the Novikov equation in $H^s$ for $s<\frac{3}{2}$ on both the line and the circle. In fact, a norm-inflation occurs (the corresponding solution has  infinite $H^s$ norm instantaneously at $t>0$, though the initial $H^s$ norm can be arbitrarily small) and gives rise to discontinuity when $\frac{5}{4}< s< \frac{3}{2},$ while for $s<\frac{5}{4}$, the solution is not unique, and either continuity or uniqueness fails when $s=\frac{5}{4}$. Very recently, Li et al. \cite{Li 2022} obtained the ill-posedness of the Novikov equation in $B_{2, \infty}^s(\mathbb{R})$ with $s>\frac{7}{2}$ due to the discontinuity.

The motivation to investigate the Novikov equation is that it can be viewed as a cubic generalization of the classical Camassa-Holm (CH) equation
\begin{eqnarray}\label{eq3}
         m_t+2u_xm=um_x, \qquad m=u-u_{xx}.
                  \end{eqnarray}
The CH equation with quadratic nonlinearity was initially derived by Fokas and Fuchssteiner \cite{Fokas 1981} as a bi-Hamiltonian system in the context of the KdV model and gained prominence after Camassa-Holm \cite{Camassa 1993} independently re-derived
it as an approximation to the Euler equations of hydrodynamics. Constantin and
Lannes \cite{Constantin 2009} later educed the CH equation from the water waves equations. The CH equation is completely integrable in the sense of having a Lax pair, a bi-Hamiltonian structure as well as possessing an infinity of conservation laws, and it also admits exact peakon solutions of the form $ce^{-|x-ct|}$ \cite{Camassa 1993, Constantin 1997, Constantin 2001, Constantin 2007, Constantin 2011,Fokas 1981}.

The local existence, uniqueness and continuity for the CH equation was initially established in the Sobolev spaces $H^s(\mathbb{R})$ with $s>\frac{3}{2}$  by Li and Olver \cite{Li 2000} as also Blanco \cite{Blanco 2001},  and then was extended by Danchin to a larger class of Besov spaces $B_{p, r}^s(\mathbb{R})$ with $s>\max\{1+\frac{1}{p}, \frac{3}{2}\}$, $1\leq p\leq \infty$  and $1\leq r< \infty$ \cite{Danchin 2001} and in the critical Besov space $B_{2, 1}^{\frac{3}{2}}(\mathbb{R})$ \cite{Danchin 2003}. It was also showed in \cite{Danchin 2003} that the CH equation is ill-posed in $B_{2, \infty}^{\frac{3}{2}}(\mathbb{R})$ due to discontinuity by using peakon solution. As the CH equation is ill-posed in the Sobolev spaces $H^s$ for $s<\frac{3}{2}$ in the sense of norm-inflation \cite{Byers 2006}, it means that the critical Sobolev exponent for well-posedness is $\frac{3}{2}$ and it was recently solved by Guo et al. in \cite{Guo 2019}, where they proved the ill-posedness of the CH equation in $B_{p, r}^{1+\frac{1}{p}}$ with $1\leq p\leq \infty$ and $1<r\leq \infty$ in the sense of norm-inflation, especially in $H^\frac{3}{2}$. Very recently, Li et al. \cite{Li 2022} showed that the CH equation is ill-posed in $B_{p, \infty}^s$ for $1\leq p \leq \infty$, $s>2+\max\{1+\frac{1}{p}, \frac{3}{2}\}$ due to the failure of continuity.

The above results do not extend clearly to the two-component Novikov system (\ref{eq1}), however, using some ideas from \cite{1Li 2022, Luo 2015}, we are able to deal with  the coupled system with these two components of the solution in different Besov spaces and establish the ill-posedness for the two-component Novikov system (\ref{eq1}) in a broader range of Besov spaces.

For studying the ill-posedness of the two-component Novikov system, it is more convenient to rewrite (\ref{eq1}) in the following equivalent nonlocal form
\begin{eqnarray}\label{eq3}
        \left\{\begin{array}{ll}
        \rho_t=u^2\rho_x+\rho uu_x,\\
         u_t=u^2u_x+\mathcal{P}(u)+\mathcal{Q}(u, \rho),\\
          \rho(0, x)=\rho_0, u(0, x)=u_0,\end{array}\right.
        \end{eqnarray}
where $\mathcal{P}(u)=\mathcal{P}_1(u)+\mathcal{P}_2(u)+\mathcal{P}_3(u),$ $\mathcal{Q}(u, \rho)=\mathcal{Q}_1(u, \rho)+\mathcal{Q}_2(u, \rho) $ and
\begin{eqnarray*}
&\;&\mathcal{P}_1(u)=\partial_x(1-\partial_x^2)^{-1}(u^3), \quad \mathcal{P}_2(u)=\frac{3}{2}\partial_x(1-\partial_x^2)^{-1}(uu_x^2), \quad \mathcal{P}_3(u)=\frac{1}{2}(1-\partial_x^2)^{-1}(u_x^3),\\
&\;&\mathcal{Q}_1(u, \rho)=-\frac{1}{2}\partial_x(1-\partial_x^2)^{-1}(u\rho^2),\quad
\mathcal{Q}_2(u, \rho)=-\frac{1}{2}(1-\partial_x^2)^{-1}(u_x\rho^2).
\end{eqnarray*}

Now our main result is stated as follows.
\begin{theorem}\label{the1.1}
Let $$s>\max\{2+\frac{1}{p},\;\; \frac{5}{2}\}, \; 1\leq p \leq \infty.$$
Then the two-component Novikov system (\ref{eq1}) is ill-posed in the Besov spaces $B_{p, \infty}^{s-1}(\mathbb{R})\times B_{p, \infty}^s(\mathbb{R})$. More precisely, there exist $(\rho_0, u_0)\in B_{p, \infty}^{s-1}\times B_{p, \infty}^s$ and a positive constant $\delta$ for which the Cauchy problem of system (\ref{eq1}) has a unique solution $(\rho, u) \in L^\infty(0, T; B_{p, \infty}^{s-1}\times B_{p, \infty}^s)$ for some $T=T(\|\rho_0\|_{B_{p, r}^{s-1} }, \|u_0\|_{B_{p, r}^s})$, while
\begin{eqnarray}\label{eq1.5}
  \liminf_{t\rightarrow 0} (\|\rho-\rho_0\|_{B_{p, \infty}^{s-1}}+ \|u-u_0\|_{B_{p, \infty}^s})\geq C\delta.
        \end{eqnarray}
 That is to say,  any energy bounded solution starting from $(\rho_0, u_0)$ does not converge back to $(\rho_0, u_0)$ in the metric of $B_{p, \infty}^{s-1}(\mathbb{R})\times B_{p, \infty}^s(\mathbb{R})$ as time goes to zero, thus results in discontinuity of the data-to-solution map and ill-posedness.
\end{theorem}

{\bf Notations}:For $\mathbf{f}=(f_1, f_2,...,f_n)\in X$,
\begin{eqnarray*}
\|\mathbf{f}\|_{X}^2=\|f_1\|_{X}^2+\|f_2\|_{X}^2+...+\|f_n\|_{X}^2.
\end{eqnarray*} 
Throughout this paper, $C$, $C_i(i=1, 2, 3, \cdot\cdot\cdot)$ stand for  universal constant which may vary from  line to line.

\section{Littlewood-Paley analysis}\label{sec2}
\setcounter{equation}{0}
In this section, we give a summary of the definition of Littlewood-Paley decomposition and nonhomogeneous Besov space, and then list some useful properties. For more details, the readers can refer to \cite{Bahouri 2011}.

There exists a couple of smooth functions $(\chi,\varphi)$ valued in $[0,1]$, such that $\chi$ is supported in the ball $\mathcal{B}\triangleq \{\xi\in\mathbb{R}^d:|\xi|\leq \frac 4 3\}$, $\varphi$ is supported in the ring $\mathcal{C}\triangleq \{\xi\in\mathbb{R}^d:\frac 3 4\leq|\xi|\leq \frac 8 3\}$. Moreover,
$$\forall\,\, \xi\in\mathbb{R}^d,\,\, \chi(\xi)+{\sum\limits_{j\geq0}\varphi(2^{-j}\xi)}=1,$$
$$\forall\,\, \xi\in\mathbb{R}^d\setminus\{0\},\,\, {\sum\limits_{j\in \mathbb{Z}}\varphi(2^{-j}\xi)}=1,$$
$$|j-j'|\geq 2\Rightarrow\textrm{Supp}\,\ \varphi(2^{-j}\cdot)\cap \textrm{Supp}\,\, \varphi(2^{-j'}\cdot)=\emptyset,$$
$$j\geq 1\Rightarrow\textrm{Supp}\,\, \chi(\cdot)\cap \textrm{Supp}\,\, \varphi(2^{-j}\cdot)=\emptyset.$$
Then, we can define the nonhomogeneous dyadic blocks $\Delta_j$ as follows:
$$\Delta_j{u}= 0,\,\, \text{if}\,\, j\leq -2,\quad
\Delta_{-1}{u}= \chi(D)u=\mathcal{F}^{-1}(\chi \mathcal{F}u),$$
$$\Delta_j{u}= \varphi(2^{-j}D)u=\mathcal{F}^{-1}(\varphi(2^{-j}\cdot)\mathcal{F}u),\,\, \text{if} \,\, j\geq 0.$$

\begin{remark}[\cite{Bahouri 2011}]\label{rem2.1}
Following the construction of $\chi$ and $\varphi$, one has $\varphi(\xi)\equiv1$ for $\frac{4}{3}\leq |\xi| \leq \frac{3}{2}.$
\end{remark}

\begin{definition}[\cite{Bahouri 2011}]\label{de2.1}
Let $s\in\mathbb{R}$ and $1\leq p,r\leq\infty$. The nonhomogeneous Besov space $B^s_{p,r}(\mathbb{R}^d)$ consists of all tempered distribution $u$ such that
\begin{align*}
||u||_{B^s_{p,r}(\mathbb{R}^d)}\triangleq \Big|\Big|(2^{js}||\Delta_j{u}||_{L^p(\mathbb{R}^d)})_{j\in \mathbb{Z}}\Big|\Big|_{\ell^r(\mathbb{Z})}<\infty.
\end{align*}
\end{definition}

Next, we list some basic lemmas and properties about Besov space which will be frequently used in proving our main result.

\begin{lemma}(\cite{Bahouri 2011})\label{lem2.1}
 (1) Algebraic properties: $\forall s>0,$ $B_{p, r}^s(\mathbb{R}^d)$ $\cap$ $L^\infty(\mathbb{R}^d)$ is a Banach algebra. $B_{p, r}^s(\mathbb{R}^d)$ is a Banach algebra $\Leftrightarrow B_{p, r}^s(\mathbb{R}^d)\hookrightarrow L^\infty(\mathbb{R}^d)\Leftrightarrow s>\frac{d}{p}$ or $s=\frac{d}{p},$ $r=1$.\\
 (2) Embedding: $B_{p, r_1}^s\hookrightarrow B_{p, r_2}^t$ for $s>t$ or $s=t $,  $r_1<r_2.$\\
(3) Let $m\in \mathbb{R}$ and $f$ be an $S^m-$ multiplier (i.e., $f: \mathbb{R}^d\rightarrow \mathbb{R}$ is smooth and satisfies that $\forall \alpha\in \mathbb{N}^d$, there exists a constant $\mathcal{C}_\alpha$ such that $|\partial^\alpha f(\xi)|\leq \mathcal{C}_\alpha(1+|\xi|)^{m-|\alpha|}$ for all $\xi \in \mathbb{R}^d$). Then the operator $f(D)$ is continuous from $B_{p, r}^s(\mathbb{R}^d)$ to $B_{p, r}^{s-m}(\mathbb{R}^d)$.\\
(4) Let  $1\leq p, r\leq \infty$ and $s>\max\{1+\frac{1}{p}, \frac{3}{2}\}$. Then  we have
$$\|uv\|_{B_{p, r}^{s-2}(\mathbb{R})}\leq C\|u\|_{B_{p, r}^{s-2}(\mathbb{R})}\|v\|_{B_{p, r}^{s-1}(\mathbb{R})}.$$
\end{lemma}

\begin{lemma}\label{lem2.2}(\cite{Bahouri 2011})
Let $s>0$, $1\leq p\leq \infty$, then we have
\begin{align*}
  \big{\|}2^{js}\|[\Delta_j, u]\partial_xv\|_{L^p}\big{\|}_{l^\infty}\leq C(\|\partial_xu\|_{L^\infty}\|v\|_{B_{p, \infty}^s}+\|\partial_xv\|_{L^\infty}\|u\|_{B_{p, \infty}^s}).
\end{align*}
Here, $[\Delta_j, u]\partial_xv=\Delta_j(u\partial_xv)-u\Delta_j(\partial_xv).$
\end{lemma}

\section{Ill-posedness }\label{sec3}
\setcounter{equation}{0}

In this section, we will give the proof of Theorem \ref{the1.1} .

Firstly, choose the initial data
\begin{align*}
  \rho_0(x)&=\sum_{n=0}^\infty2^{-n(s-1)}\phi(x)\cos(\lambda2^nx),\\
  u_0(x)&=\sum_{n=0}^\infty2^{-ns}\phi(x)\cos(\lambda2^nx),
\end{align*}
here $\lambda\in[\frac{67}{48}, \; \frac{69}{48}]$ and  $\hat{\phi}\in \mathcal{C}^\infty_0(\mathbb{R})$ is an even, real-valued and non-negative function  satisfying
\begin{numcases}{\hat{\phi}(x)=}
1, &if $|x|\leq \frac{1}{4}$,\nonumber\\
0, &if $|x|\geq \frac{1}{2}$.\nonumber
\end{numcases}

It is easy to show that
\begin{align*}
 \mathrm{ supp}\; \widehat{\phi(\cdot)\cos(\lambda2^n\cdot)}\subset \{\xi: \; -\frac{1}{2}+\lambda2^n\leq |\xi|\leq \frac{1}{2}+\lambda2^n\},
\end{align*}
then it can be verified in terms of Remark \ref{rem2.1} that for $j\geq 3$
\begin{eqnarray}\label{eq3.1}
       \Delta_j(\phi(x)\cos(\lambda2^nx))= \left\{\begin{array}{ll}
          \phi(x)\cos(\lambda2^nx),\;\; &j=n,\\
         0, \;\; &j\neq n, \end{array}\right.
        \end{eqnarray}
hence we have
\begin{align}\label{eq3.2}
  \|\rho_0\|_{B_{p, \infty}^{s-1}}\leq C, \;\;\;
 \|u_0\|_{B_{p, \infty}^s}\leq C.
\end{align}

Next, we present several estimates which play an important role in the proof of Theorem \ref{the1.1}.

\begin{lemma}\label{lem3.1}
Let $s>0$. Then for the above constructed initial data $(\rho_0, u_0)$, we have
\begin{align}
 \|u_0^2\partial_x\Delta_nu_0\|_{L^p}\geq C2^{-n(s-1)}, \label{eq3.3}\\
  \|u_0^2\partial_x\Delta_n\rho_0\|_{L^p}\geq C2^{-n(s-2)}, \label{eq3.4}
\end{align}
for   large enough $n$.
\end{lemma}
\noindent{\bf Proof} \; We just show (\ref{eq3.4}) here, since (\ref{eq3.3}) can be performed in a similar way.

 Firstly, according to (\ref{eq3.1}), we have
 \begin{align*}
   \Delta_n\rho_0=2^{-n(s-1)}\phi(x)\cos(\lambda2^nx),
 \end{align*}
hence, one has
\begin{align*}
 u_0^2\partial_x\Delta_n\rho_0=2^{-n(s-1)}u_0^2\partial_x\phi\cos(\lambda2^nx)-\lambda2^{-n(s-2)}u_0^2\phi(x)\sin(\lambda2^nx).
\end{align*}
Since $u_0^2(x)$ is a  real valued continuous function on $\mathbb{R}$, then there exists $\sigma>0$,
\begin{align}\label{eq3.5}
 |u_0^2(x)|\geq \frac{1}{2} |u_0^2(0)|=\frac{1}{2}(\sum_{n=0}^\infty2^{-ns}\phi(0))^2=\frac{2^{2s-1}\phi^2(0)}{(2^s-1)^2}, \; \mathrm{for}\; \mathrm{any}\; x\in B_\sigma(0),
\end{align}
thus we obtain from (\ref{eq3.5}) that
\begin{align*}
   \|u_0^2\partial_x\Delta_n\rho_0\|_{L^p}&\geq C2^{-n(s-2)}\|\phi(\cdot)\sin(\lambda2^n\cdot)\|_{L^p(B_\sigma(0))}-2^{-n(s-1)}\|u_0^2\partial_x\phi\cos(\lambda2^n\cdot)\|_{L^p}\\
   &\geq(C2^n-C_1)2^{-n(s-1)},
\end{align*}
which can yield (\ref{eq3.4}) by choosing $n$ large enough such that $C_1<C2^{n-1}$. Thus we finish the proof of Lemma \ref{lem3.1}.

\begin{lemma}\label{lem3.2}
Let $s>\max\{1+\frac{1}{p}, \frac{3}{2}\}$. For the above constructed initial data $(\rho_0, u_0)$, then there exists some $T=T(\|\rho_0\|_{B_{p, r}^{s-1} }, \|u_0\|_{B_{p, r}^s})$, for $0\leq t\leq T,$ we have
\begin{align}
 \|\rho(t)-\rho_0\|_{B_{p, \infty}^{s-2}}\leq Ct, \;\; \|u(t)-u_0\|_{B_{p, \infty}^{s-1}}\leq Ct.
\end{align}
\end{lemma}
\noindent{\bf Proof} \; Since $(\rho_0, u_0)\in B_{p, \infty}^{s-1}\times B_{p, \infty}^s$, according to the local existence result \cite{Luo 2015}, the Cauchy problem of system (\ref{eq1}) has a unique solution $(\rho, u) \in L^\infty(0, T; B_{p, \infty}^{s-1}\times B_{p, \infty}^s)$ for some $T=T(\|\rho_0\|_{B_{p, r}^{s-1} }, \|u_0\|_{B_{p, r}^s})$, and
\begin{align}\label{eq3.7}
 \sup_{0\leq t\leq T}(\|\rho(t)\|_{B_{p, \infty}^{s-1}}+\|u(t)\|_{B_{p, \infty}^s})\leq C(\|\rho_0\|_{B_{p, \infty}^{s-1}}+\|u_0\|_{B_{p, \infty}^s}).
\end{align}
For $t\in [0, T]$, using the differential mean value theorem, the Minkowski inequality, Lemma \ref{lem2.1} together with (\ref{eq3.7}), we obtain that
\begin{align*}
  \|\rho-\rho_0\|_{B_{p, \infty}^{s-2}}&\leq \int_0^t\|\partial_\tau\rho\|_{B_{p, \infty}^{s-2}}d\tau\\
  &\leq \int_0^t\|u^2\partial_x\rho\|_{B_{p, \infty}^{s-2}}d\tau+\int_0^t\|\rho u\partial_xu\|_{B_{p, \infty}^{s-2}}d\tau\\
   &\leq \int_0^t\|\rho\|_{B_{p, \infty}^{s-1}}\|u\|^2_{B_{p, \infty}^s}d\tau\leq Ct,
\end{align*}
and
\begin{align*}
  \|u-u_0\|_{B_{p, \infty}^{s-1}}&\leq \int_0^t\|\partial_\tau u\|_{B_{p, \infty}^{s-1}}d\tau\\
  &\leq \int_0^t\|\mathcal{P}(u)+\mathcal{Q}(u)\|_{B_{p, \infty}^{s-1}}d\tau+\int_0^t\|u^2\partial_xu\|_{B_{p, \infty}^{s-1}}d\tau\\
   &\leq \int_0^t\|\rho\|^2_{B_{p, \infty}^{s-1}}\|u\|_{B_{p, \infty}^s}d\tau+\int_0^t\|u\|^3_{B_{p, \infty}^s}d\tau\\
   &\leq Ct.
\end{align*}
Thus, we complete the proof of Lemma \ref{lem3.2}.

\begin{lemma}\label{lem3.3}
Under the assumption of Theorem \ref{the1.1}, for all $0\leq t\leq T$, we have
\begin{align}\label{eq3.8}
 \|\rho(t)-\rho_0-t\mathbf{v}_0\|_{B_{p, \infty}^{s-3}}\leq Ct^2, \;\; \|u(t)-u_0-t\mathbf{w}_0\|_{B_{p, \infty}^{s-2}}\leq Ct^2.
\end{align}
Here, $\mathbf{v}_0=u_0^2\partial_x\rho_0+\rho_0u_0\partial_xu_0,$ \; $\mathbf{w}_0=\mathcal{P}(u_0)+\mathcal{Q}(u_0)+u_0^2\partial_xu_0.$
\end{lemma}
\noindent{\bf Proof} \; For simplicity, denote
\begin{eqnarray*}
        \left\{\begin{array}{ll}
         \tilde{\rho}=\rho(t)-\rho_0-t\mathbf{v}_0,\\
          \tilde{u}=u(t)-u_0-t\mathbf{w}_0. \end{array}\right.
        \end{eqnarray*}
For $t\in [0, T]$, firstly, using the differential mean value theorem and the Minkowski inequality, we arrive at
\begin{align}
   \|\tilde{\rho}\|_{B_{p, \infty}^{s-3}} &\leq \int_0^t\|\partial_\tau\rho-\mathbf{v}_0\|_{B_{p, \infty}^{s-3}}d\tau \nonumber\\
  &\leq \int_0^t\|u^2\partial_x\rho-u_0^2\partial_x\rho_0\|_{B_{p, \infty}^{s-3}}d\tau+\int_0^t\|\rho u\partial_xu-\rho_0u_0\partial_xu_0\|_{B_{p, \infty}^{s-3}}d\tau,\label{eq3.9}
\end{align}
\begin{align}
   \|\tilde{u}\|_{B_{p, \infty}^{s-2}} &\leq \int_0^t\|\partial_\tau u-\mathbf{w}_0\|_{B_{p, \infty}^{s-2}}d\tau\nonumber\\
  &\leq \int_0^t\|\mathcal{P}(u)-\mathcal{P}(u_0)\|_{B_{p, \infty}^{s-2}}d\tau+\int_0^t\|\mathcal{Q}(u)-\mathcal{Q}(u_0)\|_{B_{p, \infty}^{s-2}}d\tau\nonumber\\
  &\;\;\;+\int_0^t\|u^2\partial_xu-u_0^2\partial_xu_0\|_{B_{p, \infty}^{s-2}}d\tau. \label{eq3.10}
\end{align}

For the term
\begin{align*}
  u^2\partial_x\rho-u_0^2\partial_x\rho_0=(u-u_0)(u+u_0)\partial_x\rho+u_0^2\partial_x(\rho-\rho_0),
\end{align*}
using (4), (2) of Lemma \ref{lem2.1}, Lemma \ref{lem3.2} and (\ref{eq3.7}), one has
\begin{align}
  \|u^2\partial_x\rho-u_0^2\partial_x\rho_0\|_{B_{p, \infty}^{s-3}}&\leq C(\|\rho\|_{B_{p. \infty}^{s-1}}\|u-u_0\|_{B_{p. \infty}^{s-1}}\|u, u_0\|_{B_{p. \infty}^s}+\|\rho-\rho_0\|_{B_{p. \infty}^{s-2}}\|u_0\|^2_{B_{p. \infty}^s})\nonumber\\
  &\leq C\tau.\label{eq3.11}
\end{align}
Similarly,
\begin{align*}
  \rho u\partial_xu-\rho_0u_0\partial_xu_0=(\rho-\rho_0)uu_x+\rho_0(u-u_0)u_x+\rho_0u_0\partial_x(u-u_0),
\end{align*}
and
\begin{align}
  \|\rho u\partial_xu-\rho_0u_0\partial_xu_0\|_{B_{p, \infty}^{s-3}}&\leq C\big(\|\rho-\rho_0\|_{B_{p. \infty}^{s-2}}\|u\|^2_{B_{p. \infty}^s}+\|u-u_0\|_{B_{p. \infty}^{s-1}}(\|u_0\|^2_{B_{p. \infty}^s}+\|\rho_0\|^2_{B_{p. \infty}^{s-1}})\big)\nonumber\\
  &\leq C\tau.\label{eq3.12}
\end{align}

With the aid of (3), (4), (2) of Lemma \ref{lem2.1}, we can find that
\begin{align}
  \|\mathcal{P}(u)-\mathcal{P}(u_0), \mathcal{Q}(u)-\mathcal{Q}(u_0)\|_{B_{p, \infty}^{s-2}}&\leq C\|u-u_0\|_{B_{p. \infty}^{s-1}}(\|u, u_0\|^2_{B_{p. \infty}^s}+\|\rho\|^2_{B_{p. \infty}^{s-1}}) \nonumber\\
 &\;\;\; +C\|\rho-\rho_0\|_{B_{p. \infty}^{s-2}}(\|\rho, \rho_0\|^2_{B_{p. \infty}^{s-1}}+\|u_0\|^2_{B_{p. \infty}^s}), \label{eq3.13}\\
 \|u^2\partial_xu-u_0^2\partial_xu_0\|_{B_{p, \infty}^{s-2}}&\leq C\|u-u_0\|_{B_{p. \infty}^{s-1}}\|u, u_0\|^2_{B_{p. \infty}^s}. \label{eq3.14}
\end{align}
Taking (\ref{eq3.11})-(\ref{eq3.12}) into (\ref{eq3.9}), (\ref{eq3.13})-(\ref{eq3.14}) into (\ref{eq3.10}) respectively, we obtain (\ref{eq3.8}). Thus we finish the proof of Lemma \ref{lem3.3}.

With Lemma \ref{lem2.2} and Lemma \ref{lem3.1}-\ref{lem3.3} at hand, we can now give the proof of Theorem \ref{the1.1}.
{\bf Proof of Theorem \ref{the1.1}} By the definition of the Besov norm, we have
\begin{align}
  \|\rho-\rho_0\|_{B_{p. \infty}^{s-1}}&\geq 2^{n(s-1)}\|\Delta_n(\rho-\rho_0)\|_{L^p}\nonumber\\
  &=2^{n(s-1)}\|\Delta_n(\tilde{\rho}+t\mathbf{v}_0)\|_{L^p}\nonumber\\
  &\geq t2^{n(s-1)}\|\Delta_n(u_0^2\partial_x\rho_0+\rho_0u_0\partial_xu_0)\|_{L^p}-2^{n(s-1)}\|\Delta_n\tilde{\rho}\|_{L^p}\nonumber\\
  &\geq t2^{n(s-1)}\|\Delta_n(u_0^2\partial_x\rho_0)\|_{L^p}-t2^{n(s-1)}\|\Delta_n(\rho_0u_0\partial_xu_0)\|_{L^p}-2^{n(s-1)}\|\Delta_n\tilde{\rho}\|_{L^p}\nonumber\\
 &\geq t2^{n(s-1)}\|\Delta_n(u_0^2\partial_x\rho_0)\|_{L^p}-Ct\|\rho_0u_0\partial_xu_0\|_{B_{p. \infty}^{s-1}}-C2^{2n}\|\tilde{\rho}\|_{B_{p, \infty}^{s-3}}.\label{eq3.15}
\end{align}
Since
\begin{align*}
 \Delta_n(u_0^2\partial_x\rho_0)&=\Delta_n(u_0^2\partial_x\rho_0)-u_0^2\partial_x\Delta_n\rho_0+u_0^2\partial_x\Delta_n\rho_0\\
 &=[\Delta_n, u_0^2\partial_x]\rho_0+u_0^2\partial_x\Delta_n\rho_0,
\end{align*}
making full use of Lemma \ref{lem2.2}, (1) of Lemma \ref{lem2.1} and (\ref{eq3.2}), we infer that
\begin{align*}
 \|2^{n(s-1)}[\Delta_n, u_0^2\partial_x]\rho_0\|_{l^\infty}&\leq C\|\partial_x(u_0^2)\|_{L^\infty}\|\rho_0\|_{B_{p, \infty}^{s-1}}+C\|\partial_x\rho_0\|_{L^\infty}\|u^2_0\|_{B_{p, \infty}^{s-1}}\\
 &\leq C\|u_0\|^2_{B_{p, \infty}^s}\|\rho_0\|_{B_{p, \infty}^{s-1}}\leq C,\\
 \|\rho_0u_0\partial_xu_0\|_{B_{p, \infty}^{s-1}}&\leq C\|\rho_0\|_{B_{p, \infty}^{s-1}}\|u_0\|^2_{B_{p, \infty}^s}\leq C.
\end{align*}
Taking above estimates into (\ref{eq3.15}), we get
\begin{align*}
  \|\rho-\rho_0\|_{B_{p, \infty}^{s-1}}\geq  t2^{n(s-1)}\|u_0^2\partial_x\Delta_n\rho_0\|_{L^p}-C_1t-C_22^{2n}\|\tilde{\rho}\|_{B_{p, \infty}^{s-3}},
\end{align*}
using Lemma \ref{lem3.1} and Lemma \ref{lem3.3} yield that
\begin{align*}
  \|\rho-\rho_0\|_{B_{p, \infty}^{s-1}}\geq  C_3t2^n-C_1t-C_22^{2n}t^2.
\end{align*}
Choosing $n$ large enough such that $C_32^n>2C_1$, then we have
\begin{align*}
  \|\rho-\rho_0\|_{B_{p, \infty}^{s-1}}\geq  \frac{C_3t2^n}{2}-C_22^{2n}t^2.
\end{align*}
As time $t$ tends to zero, picking $t2^n\approx \delta < \frac{C_3}{4C_2}$, we obtain
\begin{align*}
  \|\rho-\rho_0\|_{B_{p, \infty}^{s-1}}\geq  \frac{C_3}{2}\delta-C_2\delta^2\geq \frac{C_3}{4}\delta.
\end{align*}
 Similarly, we have
   \begin{align*}
  \|u-u_0\|_{B_{p, \infty}^s}\geq \frac{C_4}{4}\delta.
  \end{align*}
   
  This completes the proof of Theorem \ref{the1.1}.

\section*{Acknowledgments}
 This work is supported by the Natural Science Foundation of Jiangxi Province (Grant No.20212BAB211011).

\end{document}